\documentclass[12pt, a4paper]{amsart}
\usepackage{amssymb}
\usepackage{amsthm}
\usepackage{amsfonts}
\usepackage{url}
\usepackage{tikz}
\newtheorem{lemma}{\bf Lemma}[section]
\newtheorem{proposition}[lemma]{\bf Proposition}
\begin{document}
\title{Graph embeddings into Hamming spaces}
\author{Dominic van der Zypen}
\address{Swiss Armed Forces, CH-3003 Bern, Switzerland}
\email{dominic.zypen@gmail.com}
\begin{abstract}
Graph embeddings deal with injective maps from a given 
simple, undirected graph $G=(V,E)$ into a metric space, such as $\mathbb{R}^n$ with the Euclidean metric. This concept is widely studied in computer science, see \cite{ge1}, but also offers attractive research in pure graph theory \cite{ge2}. In this note
we show that any graph can be embedded into a particularly simple metric space: $\{0,1\}^n$ with the Hamming distance, for large enough $n$.
\end{abstract}
\maketitle
\section{The Hamming Graph $H(n,k)$}
We construct graph on the vertex set $\{0,1\}^n$ where $n$ is a positive integer. For $x,y \in \{0,1\}^n$ the {\it Hamming distance} of $x,y$ is the cardinality of the set $$\big\{ i \in \{0, ..., n-1\} : x(i) \neq y(i)\big\}.$$ That is, we count the positions on which $x$ and $y$ do not agree.

Fix a positive integer $k \leq n$. Two distinct elements of $\{0,1\}^n$ form an edge if their Hamming distance is at most $k$ (so they are in some sense ``close'' to each other). We denote the resulting graph on $\{0,1\}^n$ by $H(n,k)$.

We say that a finite graph $G=(V,E)$ is {\it Hamming-representable} if there are positive integers $k\leq n$ such that $G$ is isomorphic to an induced subgraph of $H(n,k)$.

As an easy example, we show that the following $3$-point graph can be embedded into $H(2,1)$:
\begin{center}
\begin{tikzpicture}
\draw (0,0.5) -- (0,2);
\draw (0.5,0) -- (2,0);
\node at (0,0) {a};
\node at (0,2.5) {b};
\node at (2.5,0) {c};
\end{tikzpicture}
\end{center}

The solution is best shown in the following picture, where it is easily seen that points connected with an edge have Hamming distance $1$ and points not connected have Hamming distance $2$:
\begin{center}
\begin{tikzpicture}
\draw (0,0.5) -- (0,2);
\draw (0.5,0) -- (2,0);
\node at (0,0) {00};
\node at (0,2.5) {01};
\node at (2.5,0) {10};
\end{tikzpicture}
\end{center}

As a further example, note that $H(n,n)$ is isomorphic to $K_{2^n}$, the complete graph on $2^n$ vertices.

Some notation: By $\mathrm{Mat}(\{0,1\}, n\times m)$ we denote the set of $n\times m$-matrices with entries in $\{0,1\}$. We identify $\mathrm{Mat}(\{0,1\}, n\times m)$ with $\{0,1\}^{nm}$ via the canonical bijection.

\section{The Result}
\begin{proposition}\label{mainprop} Every finite graph $G=(V,E)$ is Hamming-re\-pre\-sen\-ta\-ble.\end{proposition}
{\it Proof.}  We embed $G$ into $H(|E|\cdot(|V|-1),\ 2|E|-2)$. To each vertex $v$ of $G$, we will associate an $|E| \times (|V|-1)$ matrix $M_v$ with rows indexed by the edges of $G$. There will be a single $1$ in each row, with all other entries in that row equal to $0$.

If $v \in e$, then the $1$ in row $e$ of $M_v$ will be in the first column. If not, we will place a $1$ in one of the other $|V|-2$ columns, so that each of the non-endpoints of $e$ gets a $1$ in a different position of row $e$. 

If $v$ and $w$ are not joined by an edge, the Hamming distance between $M_v$ and $M_w$ is $2 |E|$ because they have no $1$'s in common; if they are joined, then the Hamming distance is $2|E|-2$.\hfill{$\Box$}.

\section{Possible use cases}
Representing graphs as subgraphs of some $H(n,k)$ can be 
useful in applications in computer science:  the
Hamming distance is computed by bitwise XOR, the
fastest operation a CPU can do. So given two
vertices represented by $n$-bit strings, it can be very quickly
determined whether they form an edge (i.e. whether their
Hamming distance is smaller than the limit given in $k$). 

Moreover, for some graphs $G=(V,E)$ with $|V|=n$ we can
represent the graph using bit strings of length ${\mathcal O}(\log n)$, making this technique potentially interesting 
for memory management.

\section{Open questions}
We define the {\em Hamming dimension} of a graph $G=(V,E)$
to be the minimum positive integer $n$ such that there is
$k\leq n$ such that $G$ can be embedded into some induced
subgraph of $H(n,k)$, and denote this by $\mathrm{dim}(G)$ Questions:
\begin{enumerate}
    \item If $G=(V,E)$ is a graph with $n=|V|$, do we necessarily have $\mathrm{dim}(G)\leq n$? If not, can we at least acheive for $\mathrm{dim}(G)$ to be ${\mathcal O}(|E|\log |V|)$?
    \item Given graphs $G,H$ what is $\mathrm{dim}(G\times H)$ in terms of $\mathrm{dim}(G), \mathrm{dim}(H)$, where $G\times H$ denotes the categorical product?
    \item How (if at all) does $\mathrm{dim}(G)$ relate to
    the chromatic number $\chi(G)$?
\end{enumerate}
\section{Acknowledgement}
I am grateful to Prof.~David Speyer of the University
of Michigan, Ann Arbor, USA, for the construction used in the proof of Proposition \ref{mainprop}.



\begin{thebibliography}{999}
\bibitem{ge1} Palash Goyal, Emilio Ferrara, {\it Graph Embedding Techniques, Applications, and Performance: A Survey}, \url{https://arxiv.org/abs/1705.02801}

\bibitem{ge2} Hongyun Cai, Vincent W. Zheng, Kevin Chen-Chuan Chang, {\it A Comprehensive Survey of Graph Embedding: Problems, Techniques and Applications}, \url{https://arxiv.org/abs/1709.07604}

\end{thebibliography}
\end{document}